\documentclass{amsart}
\usepackage{amssymb,amsmath}
\usepackage[mathscr]{euscript}
\input{epsf}

\newcounter{sec}

\newcounter{punct}[sec]

\def\punct{\refstepcounter{punct}{\arabic{sec}.\arabic{punct}.  }}

\newtheorem{theorem}{Theorem}[sec]
\newtheorem{proposition}[theorem]{Proposition}

\newtheorem{lemma}[theorem]{Lemma}

\def\COUNTERS{\addtocounter{sec}{1}
              \setcounter{punct}{0}
          \setcounter{equation}{0}
          \setcounter{theorem}{0}
          }
          
          \def\sm{\smallskip}

\begin{document}

\newcommand{\supp}{\mathop {\mathrm {supp}}\nolimits}
\newcommand{\rk}{\mathop {\mathrm {rk}}\nolimits}
\newcommand{\Aut}{\mathop {\mathrm {Aut}}\nolimits}
\newcommand{\Out}{\mathop {\mathrm {Out}}\nolimits}
\renewcommand{\Re}{\mathop {\mathrm {Re}}\nolimits}
\newcommand{\ch}{\cosh}
\newcommand{\sh}{\sinh}

\def\0{\mathbf 0}

\def\ov{\overline}
\def\wh{\widehat}
\def\wt{\widetilde}

\renewcommand{\rk}{\mathop {\mathrm {rk}}\nolimits}
\renewcommand{\Aut}{\mathop {\mathrm {Aut}}\nolimits}
\renewcommand{\Re}{\mathop {\mathrm {Re}}\nolimits}
\renewcommand{\Im}{\mathop {\mathrm {Im}}\nolimits}
\newcommand{\sgn}{\mathop {\mathrm {sgn}}\nolimits}

\def\bfa{\mathbf a}
\def\bfb{\mathbf b}
\def\bfc{\mathbf c}
\def\bfd{\mathbf d}
\def\bfe{\mathbf e}
\def\bff{\mathbf f}
\def\bfg{\mathbf g}
\def\bfh{\mathbf h}
\def\bfi{\mathbf i}
\def\bfj{\mathbf j}
\def\bfk{\mathbf k}
\def\bfl{\mathbf l}
\def\bfm{\mathbf m}
\def\bfn{\mathbf n}
\def\bfo{\mathbf o}
\def\bfp{\mathbf p}
\def\bfq{\mathbf q}
\def\bfr{\mathbf r}
\def\bfs{\mathbf s}
\def\bft{\mathbf t}
\def\bfu{\mathbf u}
\def\bfv{\mathbf v}
\def\bfw{\mathbf w}
\def\bfx{\mathbf x}
\def\bfy{\mathbf y}
\def\bfz{\mathbf z}

\def\bfA{\mathbf A}
\def\bfB{\mathbf B}
\def\bfC{\mathbf C}
\def\bfD{\mathbf D}
\def\bfE{\mathbf E}
\def\bfF{\mathbf F}
\def\bfG{\mathbf G}
\def\bfH{\mathbf H}
\def\bfI{\mathbf I}
\def\bfJ{\mathbf J}
\def\bfK{\mathbf K}
\def\bfL{\mathbf L}
\def\bfM{\mathbf M}
\def\bfN{\mathbf N}
\def\bfO{\mathbf O}
\def\bfP{\mathbf P}
\def\bfQ{\mathbf Q}
\def\bfR{\mathbf R}
\def\bfS{\mathbf S}
\def\bfT{\mathbf T}
\def\bfU{\mathbf U}
\def\bfV{\mathbf V}
\def\bfW{\mathbf W}
\def\bfX{\mathbf X}
\def\bfY{\mathbf Y}
\def\bfZ{\mathbf Z}

\def\frD{\mathfrak D}
\def\frL{\mathfrak L}
\def\frG{\mathfrak G}
\def\frg{\mathfrak g}
\def\frh{\mathfrak h}
\def\frf{\mathfrak f}
\def\frl{\mathfrak l}

\def\bfw{\mathbf w}
%%% END MATHBF
%%%%%%%%%%%%%%%%%%%%%%%%%%%%%%%
%%%%%%%%%%%%%%%%%%%%%%%%%%%%%%%%%
%%% BEGIN MATHBB

\def\R {{\mathbb R }}
 \def\C {{\mathbb C }}
  \def\Z{{\mathbb Z}}
  \def\H{{\mathbb H}}
\def\K{{\mathbb K}}
\def\N{{\mathbb N}}
\def\Q{{\mathbb Q}}
\def\A{{\mathbb A}}
\def\O {{\mathbb O }}

\def\T{\mathbb T}
\def\P{\mathbb P}

\def\G{\mathbb G}

\def\cD{\EuScript D}
\def\cL{\mathscr L}
\def\cK{\EuScript K}
\def\cM{\EuScript M}
\def\cN{\EuScript N}
\def\cP{\EuScript P}
\def\cQ{\EuScript Q}
\def\cR{\EuScript R}
\def\cW{\EuScript W}
\def\cY{\EuScript Y}
\def\cF{\EuScript F}
\def\cG{\EuScript G}
\def\cZ{\EuScript Z}
\def\cI{\EuScript I}
\def\cB{\EuScript B}
\def\cA{\EuScript A}

\def\bbA{\mathbb A}
\def\bbB{\mathbb B}
\def\bbD{\mathbb D}
\def\bbE{\mathbb E}
\def\bbF{\mathbb F}
\def\bbG{\mathbb G}
\def\bbI{\mathbb I}
\def\bbJ{\mathbb J}
\def\bbL{\mathbb L}
\def\bbM{\mathbb M}
\def\bbN{\mathbb N}
\def\bbO{\mathbb O}
\def\bbP{\mathbb P}
\def\bbQ{\mathbb Q}
\def\bbS{\mathbb S}
\def\bbT{\mathbb T}
\def\bbU{\mathbb U}
\def\bbV{\mathbb V}
\def\bbW{\mathbb W}
\def\bbX{\mathbb X}
\def\bbY{\mathbb Y}

\def\kappa{\varkappa}
\def\epsilon{\varepsilon}
\def\phi{\varphi}
\def\le{\leqslant}
\def\ge{\geqslant}

\def\B{\mathrm B}

\def\la{\langle}
\def\ra{\rangle}
\def\tri{\triangleright}

\def\lambdA{{\boldsymbol{\lambda}}}
\def\alphA{{\boldsymbol{\alpha}}}
\def\betA{{\boldsymbol{\beta}}}
\def\mU{{\boldsymbol{\mu}}}

\def\const{\mathrm{const}}
\def\rem{\mathrm{rem}}
\def\even{\mathrm{even}}
\def\SO{\mathrm{SO}}
\def\SL{\mathrm{SL}}
\def\GL{\operatorname{GL}}
\def\End{\operatorname{End}}
\def\Mor{\operatorname{Mor}}
\def\Aut{\operatorname{Aut}}
\def\inv{\operatorname{inv}}
\def\red{\operatorname{red}}
\def\Ind{\operatorname{Ind}}
\def\dom{\operatorname{dom}}
\def\im{\operatorname{im}}
\def\md{\operatorname{mod\,}}

\def\ZZ{\mathbb{Z}_{p^\mu}}

\def\cH{\EuScript{H}}
\def\cQ{\EuScript{Q}}
\def\cL{\EuScript{L}}
\def\cX{\EuScript{X}}

\def\Di{\Diamond}
\def\di{\diamond}

\def\fin{\mathrm{fin}}
\def\ThetA{\boldsymbol {\Theta}}

\def\0{\boldsymbol{0}}

\def\F{\,{\vphantom{F}}_2F_1}
\def\FF{\,{\vphantom{F}}_3F_2}
\def\H{\,\vphantom{H}^{\phantom{\star}}_2 H_2^\star}
\def\HH{\,\vphantom{H}^{\phantom{\star}}_3 H_3^\star}
\def\Ho{\,\vphantom{H}_2 H_2}

\def\disc{\mathrm{disc}}
\def\cont{\mathrm{cont}}

\def\osigma{\ov\sigma}
\def\ot{\ov t}

\begin{center}
	\bf\Large
On a $\C^2$-valued integral  index  transform
\\
and bilateral hypergeometric series

%\vspace{22pt}

\bigskip 

\sc Yury A.Neretin%
\footnote{Supported by the grant FWF,  P31591.}

\end{center}

%\vspace{22pt}

\bigskip

{\small We discuss the spectral decomposition of the hypergeometric differential operators on the line
$\Re z=1/2$, %$\mathrm{Re}\, z=1/2$
such operators arise in the problem of decomposition of tensor products of unitary representations 
of the universal covering of the group $\SL(2,\R)$. %$\mathrm{SL}\,(\mathbb R)$}
 Our main purpose is a search of natural bases in generalized eigenspaces and variants of the inversion
formula.}

\section{Introduction}

\COUNTERS

{\bf\punct The hypergeometric opertor $\cD$ on the line.}
It is well-known that classical hypergeometric systems of orthogonal polynomials
are eigenfunctions of certain differential or difference operators (see, e.g., \cite{Koe},
\cite{HTF2}, Sect. 10.8-10.13, 10.21-22, \cite{AAR}, Sect. 6.10, Ex.6.29-6.37).
On the other hand many classical integral transforms, as the Hankel transform, the Kontorovich--Lebedev transform,
the $_2F_1$-Wimp transform, the Jacobi transform (synonyms: the Olevski\u{\i} transform,
the generalized Mehler--Fock transform), etc.,
can be obtained as spectral decompositions of certain differential or difference operators
with continuous spectra (see collections of examples with differential operators in
\cite{Tit}, Chapter 4, \cite{Dun}, Sect. XIII.8).

We consider the following differential operator
 \begin{equation}
 \cD:=\frac{d}{dx} \Bigl(\frac 14+x^2\Bigr) \frac{d}{dx} +\frac{(\alpha+i\beta)^2}{4(1/2+ix)}+\frac{(\alpha-i\beta)^2}{4(1/2-ix)}+\frac 14
 \label{eq:cD}
 \end{equation}
 in $L^2(\R)$. 
 The parameters $\alpha$, $\beta$ are  real. Clearly, the replacing $(\alpha,\beta)$ by $(-\alpha,-\beta)$
 does not change the operator, so we can assume $\alpha\ge 0$. 
 
 As an algebraic expression
 $\cD$ is a hypergeometric differential operator, spectral expansions of similar operators
 produce the Jacobi polynomials (see, e.g.,\cite{Koe}, Sect.9.8) and the 'Jacobi integral transform',
 see \cite{Wey}, \cite{Tit},  Sect.4.16, \cite{Dun}, Sect. XIII.8, Theorem on p.1526, \cite{Koo1},
 \cite{Koo2}, \cite{Yak}). 
 Once more counterpart of $\cD$ was considered in \cite{MN}.

 The operator $\cD$ has a continuous spectrum
 on the half-line $\lambda\le 0$ with multiplicity two and a finite number of discrete points 
 in the domain $\lambda>0$. The explicit spectral decomposition of $\cD$ was obtained in \cite{Ner-Jacobi}.
 Since the spectrum has multiplicity, there arise a question about possible choices of natural bases in 
 spaces of solutions of the equation $\cD f=\lambda f$ for $\lambda\le 0$.
 
 \sm 

{\sc Remark.}  The operator $\cD$ appears in a natural way in the problem of decomposition of tensor products
of unitary representations of the group $\SL(2,\R)$ and its universal covering group, see \cite{Ner-Jacobi}.
Tensor products of unitary representation of $\SL(2,\R)$ were topics of many papers, in particular,
\cite{Puck}, \cite{Mol}, \cite{GKR}, \cite{Ner-Jacobi}, \cite{Gro}. However, appearance of a multiplicity
make the topic non-flexible for a further development, the same obstacles
in more serious forms arise for numerous spectral problems of non-commutative  harmonic
analysis with multiplicities. An informal purpose of the present paper is a search of an approach
to such problems.  
\hfill $\boxtimes$

\sm
 
 {\bf \punct Notation.} Denote
 $$
 \Gamma\left[\begin{matrix}
 a_1,\dots,a_n\\b_1,\dots,b_m
 \end{matrix}\right]:=\frac{\Gamma(a_1)\dots\Gamma(a_n)}{\Gamma(b_1)\dots \Gamma(b_m)}.
 $$
By
$$
\vphantom{F}_pF_q\left[\begin{matrix}
a_1,\dots,a_p\\b_1,\dots,b_q
\end{matrix};z\right]:=\sum_{n= 0}^\infty \frac{(a_1)_n\dots (a_p)_n}{(b_1)_n\dots(b_q)_n\,n!}\, z^n
$$
we denote generalized hypergeometric functions, here
$$
(a)_n:=\frac{\Gamma(a+n)}{\Gamma(a)}:=\begin{cases}
a(a+1)\dots (a+n-1)&\qquad\text{if $n\ge 0$;}
\\
\frac{1}{(a-1)\dots(a-n)}&\qquad\text{if $n\le 0$,}
\end{cases}
$$ 
is the Pochhammer symbol. By
$$
\vphantom{F}_p H_p \left[\begin{matrix}
a_1,\dots,a_p\\b_1,\dots,b_p
\end{matrix}; z\right]:=\sum_{n= -\infty}^\infty \frac{(a_1)_n\dots (a_p)_n}{(b_1)_n\dots(b_p)_n}\, z^n,
\qquad\text{where $|z|=1$.}
$$
we denote {\it bilateral hypergeometric series}, see, e.g., \cite{Sla}, Chapter 6.
If $b_p=1$, then $\sum_{n<0}$ vanishes and $(b_p)_n=(1)_n=n!$, so  we get a hypergeometric function
$\vphantom{F}_pF_{p-1}[\dots]$.
We prefer another normalization of bilateral  series
\begin{multline*}
\vphantom{F}^{\phantom{\star}}_p H_p^\star \left[\begin{matrix}
a_1,\dots,a_p\\b_1,\dots,b_p
\end{matrix}; z\right]:=
\frac{1}
{\Gamma[
1-a_1,\dots,1-a_p,\,\,b_1,\dots,b_p
]}\,\, \vphantom{F}_p H_p \left[\begin{matrix}
a_1,\dots,a_p\\b_1,\dots,b_p
\end{matrix}; z\right]=\\=
\sum_{n= -\infty}^\infty \frac{((-1)^pz)^n}{\Gamma[1-a_1-n,\dots,1-a_p-n,\,\,b_1+n,\dots, b_p+n]}.
\end{multline*}
All summands of the latter sum are well-defined (singularities are removable).
The series  absolutely converges on the circle $|z|=1$ if $\Re(\sum a_j-\sum b_j)<-1$.  If   $\Re(\sum a_j-\sum b_j)<0$, then the series conditionally
converges for $|z|=1$, $z\ne 1$.  The sum has a continuation%
\footnote{In any case,  coefficients of the  series $\vphantom{F}^{\phantom{\star}}_p H_p^\star[\dots]$ have
	polynomial growth, therefore the series always converge in the sense of distributions on the circle
	$|z|=1$.}
 to a function 
real analytic in $z=e^{i\phi}$ and meromorphic in $a_1$,\dots, $a_p$,
$b_1$, \dots, $b_p$. 

The {\it Dougall formula} (see \cite{HTF1}, formula (1.4.1), or \cite{Sla}, formula (6.1.2.1))
gives
\begin{equation}
\H\left[\begin{matrix}
a_1,a_2\\ b_1,b_2
\end{matrix};1 \right]=\frac{\Gamma[b_1+b_2-a_1-a_2-1]}{\Gamma[b_1-a_1,\,b_1-a_2,\,b_2-a_1,\,b_2-a_2]}.
\label{eq:dougall}
\end{equation}

\sm  
 
  {\bf \punct The integral transform and the inversion formula.%
  \label{ss:bases}}
 For any $\sigma\in\C$, $t\in\C$ we define a function 
 $\Phi(\sigma,t;x)=\Phi_{\alpha,\beta}(\sigma,t;x)$ on $\R$ by
 \begin{multline}
 \Phi(\sigma,t;x):=
 %\xi_{\alpha,\beta}(\sigma,t)
 \left(\tfrac12+ix\right)^{t} \left(\tfrac12-ix\right)^{-1/2-t-\sigma}
 \times\\\times
 \H\left[\begin{matrix}
 \frac{ 1-\alpha+i \beta}2+\sigma+t,\frac{1+\alpha-i \beta}2+\sigma+t\\
 1-\frac{\alpha+i\beta}{2}+t,1+\frac{\alpha+i\beta}{2}+t
 \end{matrix};-\frac{\frac12+ix}{\frac12-ix}\right],
 \label{eq:eigenfunctions}
 \end{multline} 
% where the normalizing factor%
% \footnote{A choice of the normalization is a question of a test. The present variant simplifies both
% formula ??? for $R$ and the difference operator ???.} is
% \begin{equation}
% \xi_{\alpha,\beta(\sigma,t)}:=\frac{1}%{\sin(\alpha+i\beta)\pi}
% {\Gamma\Bigl[ 1-\frac{\alpha+i\beta}2+t, 1+\frac{\alpha+i\beta}2+t,
% 	 \frac{1+\alpha-i\beta}2-t-\sigma,\frac{1-\alpha+i\beta}2-t-\sigma \Bigr]}.
% \end{equation}
 where branches of the power functions in \eqref{eq:eigenfunctions}
 are defined by the condition
 $$
 \left(\tfrac12\pm ix\right)^{\tau}\Bigr|_{x=0}=e^{\tau \ln(1/2)}.
 $$ 
 
 We have 
 $$\cD\, \Phi(\sigma,t;x)=\sigma^2 \Phi(\sigma,t;x),$$
 so for each  $\sigma$ we have a family of functions depending
 on a complex parameter $t$ in a two-dimensional space of solutions.
 For $\sigma\in i\R$ we have
  $$\Phi(\sigma,t;x)=O(|x|^{-1/2})\qquad \text{as $x\to\infty$},$$
 in this case 
 $\Phi(\sigma,t;x)$ are generalized eigenfunctions 
 of $\cD$ (see \cite{BS}, Sect. 2.2).
 
  For $f\in L^2(\R)$ we define a function
 \begin{equation*}
 J_{\alpha,\beta} f(\sigma;t):=\int_{-\infty}^\infty f(x)\,\ov{\Phi_{\alpha,\beta}(\sigma,t;x)}\,dx
 :=L^2\text{-}\!\lim_{A\to \infty} \int_{-A}^A f(x)\,\ov{\Phi_{\alpha,\beta}(\sigma,t;x)}\,dx
 \end{equation*}
 depending on $\sigma=i\nu\in i\R$, $t\in\C$. The $L^2$-$\lim$ is a limit of the family of functions
 $\phi_A:=\int_{-A}^{A}(\dots)dx$ in the sense of $L^2(\R)$.
 
   If $f\in L^2(\R)$ is compactly supported,
 then $J_{\alpha,\beta} f(\sigma;t)$ is well-defined for all 
 $\sigma$, $t\in\C^2$,  so we get  a function on
  $\C^2$ holomorphic%
  \footnote{We have $\Phi(\sigma,t;x)=O(x^{-1/2+|\Re \sigma|}\ln|x|)$ as $x\to \pm\infty$,
see \eqref{eq:PhiPsiPsi}, \eqref{eq:Psi1-as1}--\eqref{eq:Psi-as2}. The logarithm arise since formulas
\eqref{eq:Psi1-as1}--\eqref{eq:Psi-as2} are valid if $\sigma\ne 0$, 1, 2, $\dots$.
For integer $\sigma$ we come to logarithmic solutions of the hypergeometric differential
equation, see \cite{HTF1}, Subsect. 2.3.1.}
   in $\ov\sigma$, $\ov t$.
 
 We started with a function of one variable $x$ and get a function $J_{\alpha,\beta}f(\sigma,t)$ of two variables $\sigma\in i\R$,
  $t\in \C$.
 These data are overfilled,  for a reconstruction of $f$ it is sufficient to know values
 of $J_{\alpha,\beta}f(\sigma,t)$ for two values of $t$ for each $\sigma$.

\begin{theorem}
	\label{th:1}
	Let $0\le \alpha \le 1/2$.
	Consider two measurable maps $\sigma\to t(\sigma)$, $\sigma\mapsto s(\sigma)$ defined 
 for $\sigma=i\nu$, where $\nu\ge 0$. Assume that $s-t\notin\Z$ a.s. Then 
 
 \sm 
 
 {\rm a)} For 	
 $f_1\in L^2(\R)$ we have
 \begin{equation}
 f(x)=\frac1{2\pi}
 \int_0^\infty \begin{pmatrix}
 \Phi(i\nu,t(i\nu);x)& \Phi(i\nu,s(i\nu);x)
 \end{pmatrix} R
 \begin{pmatrix}
 J_{\alpha,\beta} f(i\nu,t(i\nu))\\J_{\alpha,\beta} f(i\nu,s(i\nu))
 \end{pmatrix}\,d\nu,
 \label{eq:inversion}
 \end{equation}
 where the matrix spectral density $R$ is
 \begin{multline}
 R_{t(i\nu),\,s(i\nu)}:=\\=
 \frac{\pi^4}{2\ch \pi(\beta+i\sigma)\, \ch \pi(\beta-i\sigma)\,
 	\cos\pi(\alpha-\sigma) \,\cos\pi(\alpha+\sigma)\, \Gamma[2\sigma,-2\sigma]}
 \times
 \\ \times 
 \frac{1}{\sin\pi(s-t)\,\sin\pi(\ov s-\ov t)}
 \begin{pmatrix}
 \cos\pi(\sigma+s-\ov s)&-\cos\pi(\sigma+t-\ov s)\\
 -\cos\pi(\sigma+s-\ov t) &\cos\pi(\sigma+t-\ov t)
 \end{pmatrix}.
 \label{eq:R}
 \end{multline}
 We understand the integral in \eqref{eq:inversion} as a $L^2$-limit as $B\to\infty$
 of integrals%
 \footnote{See the general statement about self-adjoint differential operators  in \cite{Dun}, Theorem XIII.5.1 and Corollary
  XIII.5.2.} $\int_0^B(\dots)d\nu$.
 
 \sm 
 
{\rm b)} For $f_1$, $f_2\in L^2(\R)$
we have the Plancherel formula
\begin{multline}
\int_{-\infty}^{\infty} f_1(x)\ov{f_2(x)}\,dx
=\\=\frac1{2\pi}\int\limits_0^{\infty} \begin{pmatrix}
J_{\alpha,\beta} f_1(i\nu,t(i\nu))&J_{\alpha,\beta} f_1(i\nu,s(i\nu))
\end{pmatrix} R_{t(i\nu),\,s(i\nu)}
\begin{pmatrix}
J_{\alpha,\beta} f_2(i\nu,t(i\nu))\\J_{\alpha,\beta} f_2(i\nu,s(i\nu))
\end{pmatrix}\,d\nu.
\label{eq:Plancherel}
\end{multline}
\end{theorem}

Theorem is proved in Section \ref{s:2}.

\sm 
 
 {\sc Remark.} It is more-or-less obvious that the matrix $R$ admits an explicit expression in the terms of
 $\Gamma$-functions. But a multiplicative structure of the matrix elements is a result of a long calculation
 and looks as a happy-end, see, e.g., transformation \eqref{eq:long} below. \hfill $\boxtimes$
 
 \sm 
 
 {\sc Remark.} Let $t=\pm\frac{\alpha+i\beta}2$. In this case
 $
 \Phi=\Phi\bigl(\sigma,\pm \tfrac{\alpha+i\beta}2;x\bigr)$
  are hypergeometric functions up 
 to simple functional factors, for this case our statement is formulated separately 
 in Proposition \ref{th:R-hypergeometric}.
 
 \sm

 {\bf\punct The case $\alpha>1/2$ and Romanovski polynomials.} This subsection contains nothing new
 comparatively \cite{Ner-Jacobi}, however it is important for understanding of our topic. 
 For $\alpha>1/2$ the operator $\cD$ has also a finite family of $L^2$-eigenfunctions
 \begin{multline*}
 \Theta^k_{\alpha,\beta}(x):=\\=
 \left(\tfrac12+ix\right)^{-(\alpha+i\beta)/2} \left(\tfrac12-ix\right)^{-(\alpha-i\beta)/2}
 \F\left[\begin{matrix}
 -k,k-2\alpha+1\\1-\alpha-i\beta
 \end{matrix};\frac12+ix \right]
=\\=
\left(\tfrac12+ix\right)^{-(\alpha+i \beta)/2}\left(\tfrac12-ix\right)^{(3\alpha+i \beta)/2-1-k}
\times\qquad\qquad\qquad\qquad\qquad\,\,\,\,
\\\times
\F\left[\begin{matrix}
k - 2 \alpha + 1, 1 - \alpha - i \beta + k\\1-\alpha-i\beta
\end{matrix};-\frac{\tfrac12+ix}{\tfrac12-ix}\right]
=\\=
\Gamma\bigl[2\alpha-k,\alpha+i\beta-k,1-\alpha-i\beta\bigr]\,
\Phi\bigl(-k+\alpha-1/2; -(\alpha+i\beta)/2;x\bigr),
\qquad
 \end{multline*}
 where $k$ ranges in integers satisfying the condition 
 \begin{equation}
 0\le k<\alpha-1/2
 \label{eq:k}
 \end{equation}
  (so for $\alpha\le 1/2$
 such functions are absent). The functions 
 $$
  R_{\alpha,\beta}^k(x):= \F\left[\begin{matrix}
 -k,k-2\alpha+1\\1-\alpha-i\beta
 \end{matrix};\frac12+ix \right]$$
  are the Romanovski polynomials \cite{Rom} (we use a nonstandard normalization),
 they are orthogonal on the line with respect to the weight 
 $$w(x)=
 (1/2+ix)^{-(\alpha+i\beta)} (1/2-ix)^{-(\alpha-i\beta)}.
 $$
 This weight decreases at infinity as a power, for this reason we have only a finite family of orthogonal
 polynomials. The  $L^2$-norms are given by (see \cite{Ask-Rom})
 \begin{multline}
 \|\Theta_{\alpha,\beta}^k\|^2_{L^2(\R)}=
\frac{1}{2\pi}\int_{-\infty}^{\infty}\frac{R_{\alpha,\beta}^k(x)\, \ov{R_{\alpha,\beta}^k(x)}\, dx}
{(1/2+ix)^{\alpha+i\beta}(1/2-ix)^{\alpha-i\beta}}
=\\=
\frac{k!\,\Gamma(2\alpha-k)}{(2n-2\alpha+1)\Gamma(\alpha+i\beta)\Gamma(\alpha-i\beta)}.
  \end{multline} 
 For $\alpha>1/2$ in the right-hand side of the inversion formula (\ref{eq:inversion}) there arise additional terms
 $$
 +\sum_k \|\Theta_{\alpha,\beta}^k\|_{L^2(\R)}^{-2}\, \la f, \Theta_{\alpha,\beta}^k\ra\,\, \Theta_{\alpha,\beta}^k(x),
 $$
where the summation is taken over $k$ satisfying  (\ref{eq:k}).
 
 \sm
 
 {\bf \punct Difference operators.}
 Next, we find the image of the operator of multiplication by $x$ under the transformation
 $J_{\alpha,\beta}$.
 
 \begin{theorem}
 	\label{th:difference}
 	Let $f(x)$ be a compactly supported integrable function on $\R$. Let 
 	the operator $J_{\alpha,\beta}$ send $f(x)$ to $F(\sigma,t)$. Then
 	$J_{\alpha,\beta}$ sends $ix\,f(x)$ to the function
 	\begin{multline}
 	\cZ F(\sigma,t)=\\=
 	\frac{(1/2+\alpha-\ov \sigma)(1/2-\alpha-\ov \sigma)(1/2+i\beta-\ov\sigma)(1/2-i\beta-\ov\sigma)}{(-2\ov\sigma)(1-2\ov\sigma)}
 	F(\sigma-1,t)+\\
 	+\frac{2i\alpha\beta}{(-1+2\ov\sigma)(1+2\ov\sigma)}F(\sigma,t)+\frac{1}{2\ov\sigma(1+2\ov\sigma)}F(\sigma+1,t).
 	\end{multline}
 \end{theorem}
 
 {\sc Remark.} Recall that in the inversion formula and in the Plancherel formula 
 the integrations are taken over the imaginary axis $\sigma\in i\R$. So $\cZ$ is a difference
 operator  in the direction transversal to the contour of integration. Similar facts take 
 place for other classical index integral transform as the Jacobi transform (see \cite{Ner-index},
 Theorem 2.1), the Kontorovich--Lebedev transform (see \cite{Ner-imaginary}, Theorem 3.2, Proposition
 3.3), the Wimp transform  (see
 \cite{Ner-imaginary}, Theorem 4.2). Moreover, Cherednik showed that the multi-dimensional
 Harish-Chandra transform sends a certain algebra of operators of multiplications to
 an algebra of difference operators
 (see \cite{Che}, \cite{vDE}).
 \hfill $\boxtimes$

\sm 

{\bf \punct The further structure of the paper.} 
Proof of Theorem \ref{th:1} is contained in Section \ref{s:2},
this section contains also two other variants of the inversion formula,
see Proposition \ref{th:R-hypergeometric}, Subsect. \ref{ss:jost}.
Theorem \ref{th:difference} is proved in Section \ref{s:3}.
The last Section \ref{s:4} contains  evaluations of  the transform of some functions.  
 
 \section{Proof of the inversion formula%
 \label{s:2}}

  \COUNTERS
  
  {\bf\punct A reduction of $\cD$ to a Schr\"odinger operator.%
  \label{ss:sch}} We consider a unitary operator
  $S:L^2(\R)\to L^2(\R)$ given by the formula
  \begin{equation}
  S f(y):=f\bigl(\tfrac12\sh y\bigr) \bigl(\tfrac12\ch y\bigr)^{1/2}.
  \label{eq:S-zamena}
  \end{equation}
It send the operator $\cD$ to the operator
$$
\cH:=\frac{d^2}{dy^2}- q(y),\qquad \text{where $q(y):=-\frac{-1+4\alpha^2-4\beta^2+8\alpha\beta\sh y}{\ch^2 y}$}
$$
(cf. \cite{Tit}, Sect. 4.16). We get a Schr\"odinger operator with a rapidly decaying potential $q(y)$,
and we can apply general statements about such operators, see, e.g., \cite{BS}, Sect. II.6.
The operator $\cH$ defined on the space $C_0^\infty(\R)$ of smooth compactly supported functions is essentially
  self-adjoint in $L^2(\R)$, see \cite{BS}, Theorem II.1.1  
  (therefore $\cD$ also is essentially self-adjoint).
The space $L^2$ splits as a direct sum $L^2(\R):=V^{\disc}\oplus V^{\cont}$ of two $\cH$-invariant subspaces
corresponding to discrete and continuous spectrum. The subspace $V^{\disc}$ is finite-dimensional,  eigenfunctions
are $L^2$-solutions of the equation $\cH f_s=s^2 f$, $s>0$ and they have asymptotics of the form
\begin{align*}
f_s(y)&=C_1 e^{-sy}\bigl(1+o(1)\bigr), \quad &\text{as $y\to+\infty$.}\\
      &=C_2 e^{sy}\bigl(1+o(1)\bigr), \quad &\text{as $y\to-\infty$.}
\end{align*}
Next, let $\sigma=i\nu\in i\R$. Consider the two-dimensional space $V_\sigma$ consisting of solutions of the equation $\cH f=\sigma^2 f$, they have asymptotics of the form
\begin{align*}
f(y)&= a\, e^{-i\nu y}\bigl(1+o(1)\bigr)+b \, e^{i\nu y}\bigl(1+o(1)\bigr),\quad &\text{as $y\to+\infty$};
\\&= c\, e^{-i\nu y}\bigl(1+o(1)\bigr)+d \, e^{i\nu y}\bigl(1+o(1)\bigr),\quad &\text{as $y\to-\infty$,}
\end{align*}
so these functions are not in $L^2$. We define an inner product in $V_\sigma$
by 
\begin{equation}
\la f_1,f_2\ra_{V_\sigma}:=\tfrac12 (a_1\ov a_2+b_1\ov b_2+c_1\ov c_2+d_1\ov d_2).
\label{eq:inner-W-sigma}
\end{equation}
Next, we define two special canonical solutions of $\cH f=\sigma^2 f$, they have asymptotics of the form
\begin{align}
\label{eq:theta1}
\theta_1(\nu;y)&= e^{i\nu y}\bigl(1+o(1)\bigr)+A(\nu)e^{-i\nu y}\bigl(1+o(1)\bigr)\quad &\text{as $y\to-\infty$;}
\\
\nonumber
   &=B(\nu)e^{i\nu y}\bigl(1+o(1)\bigr)\quad &\text{as $y\to+\infty$,}
\end{align}
and 
\begin{align}
\label{eq:theta2}
\theta_2(\nu;y)&= \qquad\qquad\qquad\quad   D(\nu)e^{-i\nu y}\bigl(1+o(1)\bigr)\quad &\text{as $y\to-\infty$;}
\\
\nonumber
&=C(\nu) e^{i\nu y}\bigl(1+o(1)\bigr)+e^{-i\nu y}\bigl(1+o(1)\bigr)\quad &\text{as $y\to+\infty$,}
\end{align}
it can be shown that the scattering matrix
$\begin{pmatrix}
A(\nu)&B(\nu)\\D(\nu)&C(\nu)
\end{pmatrix}
$
is unitary and symmetric (see \cite{BS}, Sect. II.6,  \cite{FYa}, Sect.36). For this reason,
 $\theta_1(\nu;y)$, $\theta_2(\nu;y)$ form an orthogonal basis
in $V_\sigma$ with respect to the inner product (\ref{eq:inner-W-sigma}).

Next, consider two operators, 
$$
I:L^2(\R)\to L^2(\R^+)\oplus L^2(\R^+), \qquad J:L^2(\R^+)\oplus L^2(\R^+)\to L^2(\R)
$$
given by 
$$
I:f(y)\mapsto \Bigl(\frac{1}{\sqrt{2\pi}}\int_{-\infty}^\infty f(y)\,\ov{\theta_1(\nu,y)}\,dy,
\quad \frac{1}{\sqrt{2\pi}}\int_{-\infty}^\infty f(y)\,\ov {\theta_2(\nu,y)}\,dy \Bigr )
$$
and
$$
J:\bigl(\phi_1(\nu),\phi_2(\nu)\bigr)\mapsto
 \frac{1}{\sqrt{2\pi}}\int_{-\infty}^\infty \phi_1(\nu)\,\theta_1(\nu,y)\,d\nu+
 \frac{1}{\sqrt{2\pi}}\int_{-\infty}^\infty \phi_2(\nu)\,\theta_2(\nu,y)\,d\nu.
$$
Then $\ker I=V^{\disc}$, $\im J=V^{\cont}$. The operator $I$ is a unitary operator
 $V^{\cont}\to L^2(\R)\oplus L^2(\R)$ and $J$ is the inverse operator 
 $L^2(\R)\oplus L^2(\R)\to V^{\cont}$.
 See \cite{BS}, Theorem 6.2.

We will use the statement in the following form.
{\it Let us choose {\rm (}in a measurable way{\rm )} a basis $\Psi_1(\nu;y)$, $\Psi_2(\nu;y)$ in each
$V_\sigma$. Consider the corresponding Gram matrix
\begin{equation}
\Delta(\nu):=\begin{pmatrix}
\bigl\la \Psi_1(\nu;y),\Psi_1(\nu;y)\bigr\ra& \bigl\la\Psi_1(\nu;y),\Psi_2(\nu;y)\bigr\ra\\
\bigl\la\Psi_2(\nu;y),\Psi_1(\nu;y)\bigr\ra& \bigl\la \Psi_2(\nu;y),\Psi_2(\nu;y)\bigr\ra
\vphantom{\Bigr|}
\end{pmatrix}.
\label{eq:Delta}
\end{equation}
Denote
\begin{equation}
\Xi(\nu):=\Delta(\nu)^{-1}.
\label{eq:Xi}
\end{equation}
Consider the space $C_0^\infty(\R_+)\oplus C_0^\infty(\R_+)$ equipped with the inner 
product
\begin{equation}
\Bigl\la (h_1,h_2), (h_1',h_2')\Bigr\ra:=
\frac 1{2\pi}\int_0^\infty \begin{pmatrix}
h_1(\nu)&h_2(\nu)
\end{pmatrix}  \Xi(\nu)
\begin{pmatrix}
\ov h_1^{\,\prime}(\nu)\\ \ov h_2^{\,\prime}(\nu)
\end{pmatrix}\,d\nu, 
\label{eq:XiXi}
\end{equation}
denote by $\cL[\Xi]$ the completion of $C_0^\infty(\R_+)\oplus C_0^\infty(\R_+)$
with respect to this inner product. Then the operator
$$
\cI: f\mapsto \Bigl(\bigl\la f(x), \Psi_1(\nu,x)\bigr\ra_{L^2(\R)},
 \bigl\la f(x), \Psi_2(\nu,x)\bigr\ra_{L^2(\R)}  \Bigr)
$$
is a unitary operator from $V^{\cont}$ to $\cL[\Xi]$.}

\sm 

{\bf\punct A reduction of $\cD$ to a hypergeometric differential operator.}
We set 
$$
r(x):=\left(\tfrac12+ix \right)^{(\alpha+i\beta)/2} \left(\tfrac12+ix \right)^{(\alpha-i\beta)/2} 
$$
and pass to the differential operator
$$\cB f(x):=r(x)^{-1}\cD \bigl(r(x) f(x)\bigr).$$
Next, we pass to a complex variable
$$z=\tfrac12+ix$$
and come to a new operator
$$
\cA:=-z(1-z)\frac{d^2}{dz^2}-\bigl(1+\alpha+i\beta-z(2+2\alpha)\bigr)\frac{d}{dz}+
\left(\alpha+\tfrac12\right)^2.
$$
The equation for eigenfunctions $\cA \phi=\sigma^2\phi$ becomes a special case of the hypergeometric differential equation
$$
\Bigl[z(1-z)\frac{d^2}{dz^2}+(c-z(a+b+1))\frac{d^2}{dz^2}-ab\Bigr]\phi(z)=0
$$
with
\begin{equation}
c=1+\alpha+i\beta, \qquad a=\frac12 +\alpha+\sigma,\qquad b=\frac12 +\alpha-\sigma. 
\label{eq:abc}
\end{equation}
We write two Kummer solutions 
\cite{HTF1}, (2.9.3), (2.9.20) of the hypergeometric equation
$$
(1-z)^{-a}\F\left[\begin{matrix}
a,c-b\\c
\end{matrix};\frac{z}{z-1} \right],\quad z^{1-c}(1-z)^{c-a-1}\F\left[\begin{matrix}
a+1-c,1-b\\2-c
\end{matrix};\frac{z}{z-1} \right].
$$
Substituting (\ref{eq:abc}), $z=1/2+ix$ and multiplying by $r(x)$ we get two following solutions
of the equation $\cD \psi=\sigma^2\psi$:
\begin{multline}
\Psi_1(\sigma;x)=\left(\tfrac12+ix\right)^{(\alpha+i\beta)/2} \left(\tfrac12-ix\right)^{-(\alpha+i\beta)/2-1/2-\sigma}
\times \\ \times
\F\left[\begin{matrix}
\tfrac12+\alpha+\sigma,\tfrac12+i\beta+\sigma\\1+\alpha+i\beta
\end{matrix};-\frac{1/2+ix}{1/2-ix}\right];
\label{eq:Psi1}
\end{multline} 
\begin{multline}
\Psi_2(\sigma;x)=\left(\tfrac12+ix\right)^{-(\alpha+i\beta)/2} \left(\tfrac12-ix\right)^{(\alpha+i\beta)/2-1/2-\sigma}
\times \\ \times 
\F\left[\begin{matrix}
\tfrac12-\alpha+\sigma,\tfrac12-i\beta+\sigma\\1-\alpha-i\beta
\end{matrix};-\frac{1/2+ix}{1/2-ix}\right].
\label{eq:Psi2}
\end{multline} 
These solutions are obtained one from another by a substitution $(\alpha,\beta)\leadsto(-\alpha,-\beta)$,
this substitution does not change the operator $\cD$. 

\sm

{\sc Remark.} In this place we must assume $(\alpha,\beta)\ne (0,0)$. Otherwise 
$\Psi_1$, $\Psi_2$ coincide, and we come to the logarithmic case of the hypergeometric
differential equation, see \cite{HTF1}, Sect.~2.3, \cite{AAR}, end of Section 2.3.
\hfill $\boxtimes$

\sm 

We need asymptotics of these functions as $x\to\pm\infty$. In this case the argument of hypergeometric function
tends to 1, we  apply formulas \cite{HTF1}, (2.10.1), (2.10.5),
\begin{multline*}
\F\left[\begin{matrix}
a,b\\c
\end{matrix};u\right]=\Gamma\begin{bmatrix}
c,c-a-b\\c-a,c-b
\end{bmatrix} \F\left[\begin{matrix}
a,b\\a+b-c+1
\end{matrix};1-u\right]+\\
+\Gamma\begin{bmatrix}
c,a+b-c\\a,b
\end{bmatrix} (1-u)^{c-a-b} \F\left[\begin{matrix}
c-a,c-b\\c-a-b+1
\end{matrix};1-u\right].
%\label{eq:analytic-continuation}
\end{multline*}
We have
$$
(1-u)\Bigr|_{u=-(1/2+ix)/(1/2-ix)}=\bigl(\tfrac12-ix\bigr)^{-1}.
$$
Denote 
\begin{align}
A(\alpha,\beta,\sigma)&:=\Gamma\begin{bmatrix}
1+\alpha+i\beta,\, -2\sigma\\1/2+\alpha-\sigma,\,1/2+i\beta-\sigma
\end{bmatrix};
\label{eq:A}
\\
\gamma(\alpha,\beta,\sigma)&:=\exp\Bigl\{\frac{\pi}{2}(i\alpha-\beta+i\sigma) \Bigr\}.
\label{eq:gamma}
\end{align}
Then 
\begin{multline}
\Psi_1(\sigma;x)=e^{\pi i/4}\gamma(\alpha,\beta,\sigma)\,A(\alpha,\beta,\sigma)\,x^{-1/2-\sigma}\bigl(1+o(1)\bigr) 
+\\+
e^{\pi i/4}\gamma(\alpha,\beta,-\sigma)\,A(\alpha,\beta,-\sigma)\,x^{-1/2+\sigma}\bigl(1+o(1)\bigr),
\quad \text{as $x\to +\infty$};
\label{eq:Psi1-as1}
\end{multline}
\begin{multline}
\phantom{\Psi_1(\sigma;x)}=e^{-\pi i/4}\gamma(-\alpha,-\beta,-\sigma)\,A(\alpha,\beta,\sigma)\,(-x)^{-1/2-\sigma}\bigl(1+o(1)\bigr) 
+\\+
e^{-\pi i/4}\gamma(-\alpha,-\beta,\sigma)\,A(\alpha,\beta,-\sigma)\,(-x)^{-1/2+\sigma}\bigl(1+o(1)\bigr),
\quad \text{as $x\to -\infty$}.
\label{eq:Psi-as2}
\end{multline}
For $\Psi_2(\sigma;x)$ we have a similar expression with $(\alpha,\beta)$ replaced by $(-\alpha,-\beta)$.
The formulas hold for any $\sigma\in \C$. For $\sigma\in i\R$ both summands of the asymptotic have 
the same order (and $\Psi_{1,2}(\sigma;x)$ are almost $L^2$-functions), for $\sigma\notin i\R$ one summand dominates another.

To adapt the general reasoning from Subsect. \ref{ss:sch} we must apply the unitary operator
(\ref{eq:S-zamena}), $x^{-1/2\pm\sigma}$ transform
as
\begin{align}
x^{-1/2\pm \sigma}&\leadsto \bigl(\tfrac12 \sh y\bigr)^{-1/2\pm \sigma} \bigl(\tfrac12\ch y\bigr)^{1/2}
\sim 2^{\mp\sigma} e^{\pm\sigma y}\quad \text{as $y\to+\infty$};
\label{eq:xy1}
\\
(-x)^{-1/2\pm \sigma}&\leadsto \bigl(-\tfrac12 \sh y\bigr)^{-1/2\pm \sigma} \bigl(\tfrac12\ch y\bigr)^{1/2}
\sim 2^{\mp\sigma} e^{\mp\sigma y}\quad \text{as $y\to-\infty$}.
\label{eq:xy2}
\end{align}

{\bf \punct The Gram matrix for the hypergeometric eigenfunctions.} Let $\sigma\in i\R$. 
Our next purpose is to evaluate the matrices $\Delta$ and $\Xi=\Delta^{-1}$ 
(see \eqref{eq:Delta}) for the eigenfunctions $\Psi_1$, $\Psi_2$ given by \eqref{eq:Psi1}, 
 \eqref{eq:Psi2}.
 
 \begin{lemma}
 	\label{l:Delta}
 	{\rm a)} The matrix elements of the Gram matrix $\Delta$ for eigenfunctions $\Psi_1$, $\Psi_2$ are
 	\begin{align*}
 	\Delta_{11}&=\frac 2\pi \ch \pi(\beta-i\sigma)\,\ch \pi(\beta+i\sigma) \,
 	\Gamma\begin{bmatrix}
 	1+\alpha+i\beta,\, 1+\alpha-i\beta,\,2\sigma,-2\sigma\\
 	1/2+\alpha-\sigma,\, 1/2+\alpha+\sigma
 	\end{bmatrix};
 	\\
 	\Delta_{12}&=\frac 2\pi \cos \pi(\alpha-\sigma)\,\cos \pi(\alpha+\sigma)\,
 	\Gamma\begin{bmatrix}
 	1-\alpha+i\beta,\, 1+\alpha+i\beta,\, 2\sigma,\,-2\sigma\\
 	1/2+i\beta-\sigma,\,1/2+i\beta+\sigma
 	\end{bmatrix};\\
 \Delta_{21}&=\frac 2\pi \cos \pi(\alpha-\sigma)\,\cos \pi(\alpha+\sigma)\,
 \Gamma\begin{bmatrix}
 1-\alpha-i\beta,\,  1+\alpha-i\beta,\, 2\sigma,\,-2\sigma\\
 1/2-i\beta-\sigma,\,1/2-i\beta+\sigma
 \end{bmatrix};\\	
 \Delta_{22}&=\frac 2\pi \ch \pi(\beta-i\sigma)\,\ch \pi(\beta+i\sigma) \,
 \Gamma\begin{bmatrix}
 1-\alpha-i\beta,\,1-\alpha+i\beta,\, 2\sigma,\,-2\sigma\\
 1/2-\alpha-\sigma,\,1/2-\alpha+\sigma
 \end{bmatrix}.
 	\end{align*}
 	
 {\rm b)} The determinant of $\Delta$ is
 \begin{multline*}
 \det \Delta=\frac 4{\pi^2}
  \cos\pi(\alpha - \sigma)\,
  \cos \pi (\alpha +  \sigma)\, \ch \pi (\beta - i\sigma)\, \ch \pi(\beta + i\sigma)
  \times\\ \times
  (\alpha^2 +  \beta^2)\,	\Gamma[2\sigma,-2\sigma]^2.
  \end{multline*}
 \end{lemma}

{\sc Proof.}
a) We present a calculation of $\Delta_{11}$,
\begin{multline}
\Delta_{11}=\frac12\bigl(|\gamma(\alpha,\beta,\sigma)|^2+|\gamma(-\alpha,-\beta,-\sigma)|^2\bigr)
\,\,|A(\alpha,\beta,\sigma)|^2+\\+
\frac12\bigl(|\gamma(\alpha,\beta,-\sigma)|^2+|\gamma(-\alpha,-\beta,+\sigma)|^2\bigr)
\,\,|A(\alpha,\beta,-\sigma)|^2.
\label{eq:Delta11}
\end{multline}
Let us evaluate the first summand. We have
$$
\tfrac12\bigl(|\gamma(\alpha,\beta,\sigma)|^2+|\gamma(-\alpha,-\beta,-\sigma)|^2\bigr)=
\tfrac12 (e^{\pi(-\beta+i\sigma)}+e^{-\pi(-\beta+i\sigma)})=\ch\pi(\beta-i\sigma)
$$
and
\begin{multline*}
|A(\alpha,\beta,\sigma)|^2=\Gamma\begin{bmatrix}
1+\alpha+i\beta,\, -2\sigma\\1/2+\alpha-\sigma,\,1/2+i\beta-\sigma
\end{bmatrix}
\Gamma\begin{bmatrix}
1+\alpha-i\beta,\, 2\sigma\\1/2+\alpha+\sigma,\,1/2-i\beta+\sigma
\end{bmatrix}=\\=
\frac1\pi\ch\pi(\beta+i\sigma)\,\Gamma\begin{bmatrix}
1+\alpha+i\beta,\,1+\alpha-i\beta,\, -2\sigma,\,2\sigma\\1/2+\alpha-\sigma,\,1/2+\alpha+\sigma
\end{bmatrix}.
\end{multline*}
We see that the first summand in \eqref{eq:Delta11} is even in $\sigma$. Therefore it is equal
to the second summand, and we come to the final expression.

Evaluations of other matrix elements are similar.

\sm  

b) Evaluating $\det\Delta=\Delta_{11}\Delta_{22}-\Delta_{12}\Delta_{21}$ we meet
the following subexpressions:
\begin{align*}
&\frac{1}{\Gamma[\tfrac12+\alpha+\sigma,\, \tfrac12+\alpha-\sigma,\,\tfrac12-\alpha+\sigma,\,\tfrac12-\alpha-\sigma]}
=\frac{1}{\pi^2} \cos\pi(\alpha-\sigma)\,\cos\pi(\alpha+\sigma);
\\
&\frac{1}{\Gamma[\tfrac12+i\beta+\sigma,\,\tfrac12+i\beta-\sigma,\,\tfrac12-i\beta+\sigma,\,\tfrac12-i\beta-\sigma]}
=\\
&\qquad\qquad\qquad\qquad\qquad\qquad\qquad\qquad=\frac{1}{\pi^2}\ch\pi(\beta-i\sigma)\,\ch\pi(\beta+i\sigma);
\\
&\Gamma[1+\alpha+i\beta,\,1+\alpha-i\beta,\,1-\alpha+i\beta,\,1-\alpha-i\beta]=
\frac{\pi^2(\alpha+i\beta)(\alpha-i\beta)}{\sin\pi(\alpha+i\beta)\,\sin\pi(\alpha-i\beta)}.
\end{align*}
Applying these transformations we get the following expression for $\det\Delta$:
\begin{multline*}
\!\!\!\!\!\!\frac{4}{\pi^2}
\frac{(\alpha^2+\beta^2)\cos\pi(\alpha-\sigma)\cos\pi(\alpha+\sigma)
	\ch\pi(\beta-i\sigma)\ch\pi(\beta+i\sigma)
\Gamma[2\sigma,-2\sigma]^2}
{\sin\pi(\alpha+i\beta)\,\sin\pi(\alpha-i\beta)}
\times\\\times
\Bigl\{\ch\pi(\beta-i\sigma)\ch\pi(\beta+i\sigma)-\cos\pi(\alpha-\sigma)\cos\pi(\alpha+\sigma) \Bigr\}.
\end{multline*}
Simplifying the expression in the curly brackets we get
$$
\Bigl\{\dots  \Bigr \}=\sin\pi(\alpha+i\beta)\,\sin\pi(\alpha-i\beta)
$$
and we come to the final expression.
\hfill$\square$

\sm

Now we write the matrix $\Delta^{-1}$ in a straightforward way and get the following statement, see \cite{Ner-Jacobi}:

\begin{proposition}
	\label{th:R-hypergeometric} Let $(\alpha,\beta)\ne (0,0)$.
	Then for the eigenfunctions $\Psi_1(x)$, $\Psi_2(x)$ given by \eqref{eq:Psi1}--\eqref{eq:Psi2}
	the spectral matrix $\Xi$ in \eqref{eq:XiXi} is given by%
	\footnote{We write each matrix element as a two-line formula, this allows to obtain
	a readable expression.}
$$
\frac{1}{2\pi\Gamma[2\sigma,-2\sigma]}
\begin{pmatrix}
\begin{matrix}
\Gamma[\frac12+\alpha+\sigma,\frac12+\alpha-\sigma]\times\,\,\\
\,\,\,\,\,\,\times\Gamma[-\alpha-i\beta,-\alpha+i\beta]
\end{matrix}_{\vphantom{\Bigl|}}&
\begin{matrix}
\Gamma[\frac12-i\beta+\sigma,\frac12-i\beta-\sigma]\times\,\\
\,\,\,\,\,\,\times\Gamma[-\alpha+i\beta,\alpha+i\beta]
\end{matrix}\\
\begin{matrix}
\Gamma[\frac12+i\beta+\sigma,\frac12+i\beta-\sigma]\times\,\\
\,\,\,\,\,\,\times\Gamma[-\alpha-i\beta,\alpha-i\beta]
\end{matrix}^{\vphantom{A^A}}&
\begin{matrix}
\Gamma[\frac12-\alpha+\sigma,\frac12-\alpha-\sigma]\times\,\,\\
\,\,\,\,\,\,\times\Gamma[\alpha-i\beta,\alpha+i\beta]
\end{matrix}
\end{pmatrix}.
$$	
\end{proposition}

{\bf \punct Bilateral hypergeometric functions $\H$.}
It is easy to see that functions $\Ho\left[\begin{matrix}
a_1,a_2\\b_1,b_2
\end{matrix};z\right]$
satisfy the differential equation
$$
\Bigl\{ z\Bigl(z\frac{d}{dz}+a_1 \Bigr)\Bigl(z\frac{d}{dz}+a_2 \Bigr)-
\Bigl(z\frac{d}{dz}+b_1-1 \Bigr) \Bigl(z\frac{d}{dz}+b_2-1\Bigr)   \Bigr\} \,F(z)=0
$$
(cf. \cite{Sla}, formula (2.1.2.1)). Functions $\H$ differ from $\Ho$ by constant factors.
 Moreover for any $t\in\C$
the function
$$J_t(z):=(-z)^{t}\H\left[\begin{matrix}
a_1+t,a_2+t\\b_1+t,b_2+t
\end{matrix};z\right]
$$
satisfy the same differential equation (we assume that $(-z)^t\bigr|_{z=-1}=1$). Therefore 
any 3 functions $J_{t_1}(z)$, $J_{t_2}(z)$, $J_{t_3}(z)$  are linear dependent, i.e.,
\begin{equation}
C_1 J_{t_1}(z) +
C_2 J_{t_2}(z)+C_3J_{t_3}(z)=0
\label{eq:C1C2C3}
\end{equation}
for some $C_1$, $C_2$, $C_3$. In fact,
see \cite{Aga},
\begin{equation}
\sin\pi(t_2-t_3)\,J_{t_1}(z) +
\sin\pi(t_3-t_1)\,J_{t_2}(z)+\sin\pi(t_1-t_2)\,J_{t_3}(z)=0.
\label{eq:2H2-dependence}
\end{equation}

{\sc Remark.}
These coefficients $C_1$, $C_2$, $C_3$ of the linear dependence
 can be derived in the following way.
The Dougall formula \eqref{eq:dougall}
provides us an explicit value for any $\H(z)$ at $z=1$. We substitute $z=e^{0_+i}$ and 
$z=e^{2\pi_- i}$ and get two equations for the coefficients.
\hfill $\boxtimes$

\sm 

Setting
$$
t_1=0,\qquad t_2=1-b_1,\quad t_3=1-b_2
$$
to \eqref{eq:2H2-dependence} we get an expression of an arbitrary function $\H(z)$ in terms
of Gauss hypergeometric functions.

\sm 

In particular, we get an expression for the functions
\begin{multline*}
\Phi(\sigma,t;x):=
\left(\tfrac12+ix\right)^{t} \left(\tfrac12-ix\right)^{-1/2-t-\sigma}
\times\\\times
\H\left[\begin{matrix}
\frac{ 1-\alpha+i \beta}2+\sigma+t,\frac{1+\alpha-i \beta}2+\sigma+t\\
1-\frac{\alpha+i\beta}{2}+t,1+\frac{\alpha+i\beta}{2}+t
\end{matrix};-\frac{\frac12+ix}{\frac12-ix}\right],
\end{multline*}
 defined in Subsect. \ref{ss:bases}.
 Namely,
 \begin{multline}
 \label{eq:PhiPsiPsi}
 \Phi(\sigma,t;x)\sin\pi(\alpha+i\beta)  =\\= C_1(\sigma) \Psi_1(\sigma;x)\sin\pi\bigl(t+\tfrac{\alpha+i\beta}2 \bigr)  
 +C_2(\sigma)\Psi_2(\sigma;x)\sin\pi\bigl(-t+\tfrac{\alpha+i\beta}2 \bigr), 
 \end{multline}
where
\begin{align}
C_1(\sigma)&=\frac 1{\Gamma[\frac12-i\beta-\sigma,\frac12-\alpha-\sigma,1+\alpha+i\beta]}:=C(\alpha,\beta,\sigma);
\label{eq:C1}
\\
C_2(\sigma)&=\frac 1{\Gamma[\frac12+i\beta-\sigma,\frac12+\alpha-\sigma,1-\alpha-i\beta]}
=C(-\alpha,-\beta,\sigma).
\label{eq:C2}
\end{align}

\begin{lemma}
	\begin{equation}
	\bigl\la \Phi(\sigma,t;x), \Phi(\sigma,s;x)\bigr\ra_{V_\sigma}=
	M\cdot \cos \pi(\sigma+t-\ov s), 
	\label{eq:PhiPhi}
	\end{equation}
	where
$$
M=\ch\pi(\beta+i\sigma)\,\ch\pi(\beta-i\sigma)\,\cos\pi(\alpha+\sigma) \,\cos\pi(\alpha-\sigma)\,
\Gamma[2\sigma,-2\sigma].
$$	
\end{lemma}

{\sc Proof.} Let $\Delta$ be the Gram matrix of the eigenfunctions $\Psi_1$, $\Psi_2$, see Lemma \ref{l:Delta}.  
Then
\begin{equation*}
\begin{pmatrix}
C_1 \Delta_{11} \ov C_1&C_1 \Delta_{12} \ov C_2\\
C_2 \Delta_{21} \ov C_1&C_2 \Delta_{22} \ov C_2
\end{pmatrix}
=\frac{2}{\pi^4}M%\cdot\sin\pi(\alpha+i\beta)\,\sin\pi(\alpha-i\beta)
\cdot S,
\end{equation*}
where
$$
S=
\begin{pmatrix}
\ch\pi(\beta-i\sigma)&\cos\pi(\alpha+\sigma)\\
\cos\pi(\alpha-\sigma)&\ch\pi(\beta+i\sigma)
\end{pmatrix}.
$$
Let us verify the identity for the first matrix element:
\begin{multline*}
C_1\Delta_{11}\ov C_1=
\frac 2\pi \ch \pi(\beta-i\sigma)\,\ch \pi(\beta+i\sigma) \,
\Gamma\begin{bmatrix}
1+\alpha+i\beta,\, 1+\alpha-i\beta,\,2\sigma,-2\sigma\\
1/2+\alpha-\sigma,\, 1/2+\alpha+\sigma
\end{bmatrix}
\times\\ \times 
\frac 1{\Gamma[\frac12-i\beta-\sigma,\frac12-\alpha-\sigma,1+\alpha+i\beta]}
\cdot \frac 1{\Gamma[\frac12+i\beta+\sigma,\frac12-\alpha+\sigma,1+\alpha-i\beta]}.
\end{multline*}
The product of three $\Gamma$-factors is
$$
\frac{1}{\pi^3} \cos\pi(\alpha+\sigma)\,\cos\pi(\alpha-\sigma)\,\Gamma[-2\sigma,2\sigma]\, \ch\pi(\beta-i\sigma),
$$
and we come to the desired expression.

 Now we are ready to evaluate
\begin{multline}
\label{eq:PhiPhi1}
\bigl\la \Phi(\sigma,t;x), \Phi(\sigma,s;x)\bigr\ra_{V_\sigma}=\frac{2}{\pi^4}
\frac{M}{\sin\pi(\alpha+i\beta)\,\sin\pi(\alpha-i\beta)}
\times \\ \times
\left\{
\begin{pmatrix}
\sin\pi(\frac{\alpha+i\beta}{2}+t)&\sin\pi(\frac{\alpha+i\beta}{2}-t)
\end{pmatrix}\,S\,
\begin{pmatrix}
\sin\pi(\frac{\alpha-i\beta}{2}+\ov s)\\ \sin\pi(\frac{\alpha-i\beta}{2}-\ov s)
\end{pmatrix}
 \right\}.
\end{multline}
The expression in the curly bracket is
\begin{multline}
\label{eq:long}
\Bigl\{\dotsc\Bigr \}=\sin\pi(\tfrac{\alpha+i\beta}{2}+t)\,\ch\pi(\beta-i\sigma) \,\sin\pi(\tfrac{\alpha-i\beta}{2}+\ov s)
+\\+
\sin\pi(\tfrac{\alpha+i\beta}{2}+t)\,\cos\pi(\alpha+\sigma)\,\sin\pi(\tfrac{\alpha-i\beta}{2}-\ov s)
+\\+
\sin\pi(\tfrac{\alpha+i\beta}{2}-t)\,\cos\pi(\alpha-\sigma)\,\sin\pi(\tfrac{\alpha-i\beta}{2}+\ov s)
+\\+
\sin\pi(\tfrac{\alpha+i\beta}{2}-t)\,\ch\pi(\beta+i\sigma)\,\sin\pi(\tfrac{\alpha-i\beta}{2}-\ov s)
=\\=
\cos\pi(\sigma+t-\ov s)\,\sin\pi(\alpha+i\beta)\,\sin\pi(\alpha-i\beta),
\end{multline}
this implies the statement of the lemma. The last identity is not obvious, but when written,
it admits a straightforward verification.
\hfill $\square$

\sm

{\sc Proof of Theorem \ref{th:1}.} Thus, the Gram matrix
of $\Phi(\sigma,t;x)$ and $\Phi(\sigma,s;x)$
is
$$
\frac{2}{\pi^4}\cdot M \begin{pmatrix}
\cos\pi(\sigma+t-\ov t)&\cos\pi(\sigma+t-\ov s)\\
\cos\pi(\sigma+s-\ov t)&\cos\pi(\sigma+s-\ov s)
\end{pmatrix}.
$$
The inverse matrix is 
$$
\frac{\pi^4}{2} \frac1M\cdot \frac{1}{\sin\pi(s-t)\,\sin\pi(\ov s-\ov t)}
\begin{pmatrix}
\cos\pi(\sigma+s-\ov s)&-\cos\pi(\sigma+t-\ov s)\\
-\cos\pi(\sigma+s-\ov t)&\cos\pi(\sigma+t-\ov t)
\end{pmatrix}
$$
and we come to the formula \eqref{eq:Plancherel}.
\hfill $\square$

\sm

{\bf\punct A generalized orthogonal system.%
\label{ss:jost}} For completeness we present formulas for 
the eigenfunctions $\theta_1$, $\theta_2$, see \eqref{eq:theta1}--\eqref{eq:theta2}.
Set
\begin{multline*}
\theta_1(\sigma;x):=\\=
-\frac{e^{\frac{3\pi i}4}}{2\pi}\,
\Bigl(\mu(\alpha,\beta,\sigma)M(\alpha,\beta,\sigma) \Psi_1(\sigma;x)+ \mu(-\alpha,-\beta,\sigma)M(-\alpha,-\beta,\sigma) \Psi_2(\sigma;x)\Bigr);
\end{multline*}
\begin{multline*}
\theta_2(\sigma;x):=\\
=\frac{e^{\frac{\pi i}4}}{2\pi}\,
\Bigl(\mu(-\alpha,-\beta,-\sigma)M(\alpha,\beta,\sigma) \Psi_1(\sigma;x)+ \mu(\alpha,\beta,-\sigma)M(-\alpha,-\beta,\sigma) \Psi_2(\sigma;x)\Bigr),
\end{multline*}
where
\begin{align*}
\mu(\alpha,\beta,\sigma)&:=e^{\frac{\pi}2(-i\alpha+\beta+i\sigma)};
\\
M(\alpha,\beta,\sigma)&:=
\Gamma\left[\begin{matrix}
-\alpha-i\beta,\frac12 + \alpha - \sigma,1/2 + i \beta - \sigma\\
-2\sigma
\end{matrix} \right].
\end{align*}
Then $\theta_1$, $\theta_2$ have the following asymptotics at infinity 
(see \eqref{eq:xy1}--\eqref{eq:xy2}):
\begin{align*}
\theta_1(\sigma;x)&= (-x)^{-\frac12-\sigma}(1+O(x^{-1}))+A(\sigma)(-x)^{-\frac12+\sigma}(1+O(x^{-1}))
\,\, &\text{as $x\to-\infty$;}
\\
\nonumber
&=B(\sigma)x^{-\frac12+\sigma}(1+O(x^{-1}))\,\, &\text{as $x\to+\infty$,}
\end{align*}
and 
\begin{align*}
\theta_2(\sigma;y)&= \qquad\qquad\qquad\quad   D(\sigma)(-x)^{-\frac12+\sigma}(1+O(x^{-1}))\,
\, &\text{as $x\to-\infty$;}
\\
\nonumber
&=C(\sigma) x^{-\frac12+\sigma} (1+O(x^{-1}))+x^{-\frac12-\sigma}(1+O(x^{-1}))\,\, &\text{as $x\to+\infty$,}
\end{align*}
where the elements of the scattering matrix are given by
\begin{align*}
A&:=\frac{1}{2\pi^2}\bigl(e^{-\pi\beta}\cos\pi(\alpha - \sigma)+e^{\pi\beta}\cos\pi(\alpha + \sigma) \bigr)
\Gamma\left[\begin{matrix}
2\sigma\\-2\sigma
\end{matrix} \right]
\cdot\cG;
\\
B&=D:=\frac{1}{2\pi \Gamma[1 - 2 \sigma,-2 \sigma]} \cdot \cG;
\\
C&:= \frac{1}{2\pi^2}\bigl(e^{\pi\beta}\cos\pi(\alpha - \sigma)+e^{-\pi\beta}\cos\pi(\alpha + \sigma) \bigr)
\Gamma\left[\begin{matrix}
2\sigma\\-2\sigma
\end{matrix} \right]\cdot \cG,
\end{align*}
and
$$
\cG:=\Gamma\bigl[\tfrac12 - \alpha - \sigma,\, \tfrac12 + \alpha - \sigma,\,\tfrac12 - i\beta - \sigma,\,
\tfrac12 + i\beta - \sigma\bigr].
$$

For such functions $\theta_1$, $\theta_2$ the matrix $\Xi$ in \eqref{eq:XiXi}
is $\begin{pmatrix}
1&0\\0&1
\end{pmatrix}$,
but we pay for this by longer and less flexible expressions for eigenfunctions.

\sm

{\bf \punct The case $\alpha=0$, $\beta=0$.} For this case the calculations of this section are not
valid, but we can easily apply continuity arguments. Our  final formula \eqref{eq:R} follows from \eqref{eq:PhiPhi}.
To extend the latter formula to our case, it is sufficient to show, that coefficients  at high terms of asymptotics $|x|^{-1/2\pm\sigma}$ of $\Phi_{\alpha,\beta}(\sigma,t;x)$ at infinities are continuous at the point $(\alpha,\beta)=(0,0)$
for $\sigma=i\nu$, where $\nu>0$. These coefficients can be easily written explicitly with formulas
\eqref{eq:PhiPsiPsi}, \eqref{eq:Psi1-as1}-\eqref{eq:Psi-as2}. For instance, the coefficient in front
of $x^{-1/2-\sigma}$ as $x\to+\infty$ 
is
\begin{multline}
\label{eq:coeff}
\frac{1}{\sin \pi(\alpha+i\beta)/2}
\Bigl\{\gamma(\alpha,\beta,\sigma)\, A(\alpha,\beta,\sigma)\, C(\alpha,\beta,\sigma) 
\sin\pi\bigl( t+\tfrac{\alpha+i\beta}{2}\bigr) -
\\-\gamma(-\alpha,-\beta,\sigma)\, A(-\alpha,-\beta,\sigma)\, C(-\alpha,-\beta,\sigma) 
\sin \pi\bigl( t-\tfrac{\alpha+i\beta}{2}\bigr)  \Bigr\},
\end{multline}
where $\gamma(\dots)$, $A(\dots)$, $C(\dots)$ are defined by formulas
\eqref{eq:A}-\eqref{eq:gamma},
\eqref{eq:C1}.  Substitute $\alpha=-i\beta$ to the bracket $\Bigl\{\dots\Bigr \}$.
It is easy to see that
$$
\gamma(\alpha,\beta,\sigma)\Bigr|_{\alpha=-i\beta}=\gamma(-\alpha,-\beta,\sigma)\Bigr|_{\alpha=-i\beta},
$$
similar identities take place also for $A(\dots)$, $C(\dots)$.
Therefore  for $\alpha=-i\beta$ the expression $\Bigl\{\dots\Bigr \}$ is zero.
So the singularity on the surface $\alpha+i\beta=0$ in \eqref{eq:coeff} is removable
and the whole expression is continuous.

\section{The difference operator%
\label{s:3}}

\COUNTERS

The topic of this section is the proof of Theorem \ref{th:difference}. In fact, we must show that 
the kernel $F=\Phi(\sigma,t;x)$ satisfies the equation
	\begin{multline}
	\label{eq:difference-F}
-ix\,F(\sigma,t;x)=\\=
\frac{(1/2+\alpha- \sigma)(1/2-\alpha- \sigma)(1/2+i\beta-\sigma)(1/2-i\beta-\sigma)}{(-2\sigma)(1-2\sigma)}
F(\sigma-1,t;x)
-\\
-\frac{2i\alpha\beta}{(-1+2\sigma)(1+2\sigma)}F(\sigma,t;x)
+\frac{1}{2\sigma(1+2\sigma)}F(\sigma+1,t;x)
\end{multline}
(the variable $t$ is absent in the coefficients).

We use the expression \eqref{eq:PhiPsiPsi} for $\Phi$, it is sufficient to show that two terms
$C_1\Psi_1(x)$, $C_2\Psi_2(x)$ satisfy the same difference equation.
We write $\Psi_1(x)$ as
\begin{equation*}
\Psi_1(\sigma;x)=
\left(\tfrac12+ix \right)^{(\alpha+i\beta)/2} \left(\tfrac12-ix \right)^{(\alpha-i\beta)/2}
\F\left[\begin{matrix}
\frac12+\alpha+\sigma,\frac12+\alpha-\sigma\\1+\alpha+ i\beta 
\end{matrix};\frac12+ix \right],
\end{equation*}
see \cite{HTF1}, formulas (2.9.1), (2.9.3). The expression for $\Psi_2(\sigma;x)$
is obtained by replacing $(\alpha,\beta)\leadsto(-\alpha,-\beta)$

By \cite{Ner-index}, formula (2.3), the Gauss hypergeometric function satisfies the following 
contiguous relation
\begin{multline*}
-y\F\left[\begin{matrix}
p,q\\r
\end{matrix};y \right]= \frac{q(r-p)}{(q-p)(1+q-p)}\F\left[\begin{matrix}
p-1,q+1\\r
\end{matrix};y \right]
-\\-
\Bigl(\frac{q(r-p)}{(q-p)(1+q-p)}+\frac{p(r-q)}{(p-q)(1+p-q)}\Bigr)\,\F\left[\begin{matrix}
p,q\\r
\end{matrix};y\right]
+\\
+\frac{p(r-q)}{(p-q)(1+p-q)}\F\left[\begin{matrix}
p+1,q-1\\r
\end{matrix};y \right].
\end{multline*}
Therefore $G(\sigma;x)=\Psi_1(\sigma;x)$ satisfies the difference equation
\begin{multline}
\label{eq:difference-G}
-\left(\tfrac12+ix \right) G(\sigma,x)=
\frac{(\tfrac12 +\alpha-\sigma)(\tfrac12 +i\beta-\sigma)}{-2\sigma(1-2\sigma)}G(\sigma-1,x)
-\\-
\left\{\frac{(\tfrac12 +\alpha-\sigma)(\tfrac12 +i\beta-\sigma)}{-2\sigma(1-2\sigma)}+
\frac{(\tfrac12 +\alpha+\sigma)(\tfrac12 +i\beta+\sigma)}{2\sigma(1+2\sigma)} \right\}
G(\sigma,x)
+\\+
\frac{(\tfrac12 +\alpha+\sigma)(\tfrac12 +i\beta+\sigma)}{2\sigma(1+2\sigma)}F(\sigma+1,x)
\end{multline}
The expression in the curly brackets can be transformed as
$$
\Bigl\{\dots\Bigr \}=
\frac12+ \frac{2i\alpha\beta}{(-1+2\sigma)(1+2\sigma)}.
$$
If $G(\sigma,x)$ satisfies the difference equation \eqref{eq:difference-G},
then $F(\sigma;x)=C_1(\sigma)G(\sigma;x)$
satisfies equation \eqref{eq:difference-F}. So $C_1(\sigma)\Psi_1(\sigma;x)$ satisfies
\eqref{eq:difference-F}. Since expression \eqref{eq:difference-F}
is invariant with respect to the transformation
$(\alpha,\beta)\leadsto(-\alpha,-\beta)$, the summand $C_2(\sigma)\Psi_2(\sigma;x)$
also satisfies the difference equation.

\section{Some evaluations%
\label{s:4}}

\COUNTERS

There are many explicit evaluations for the Jacobi transform, this allows to use it as a tool
for obtaining non-trivial properties of special functions, see \cite{Koo2}, \cite{Ner-beta},
\cite{Ner-add} (see also, \cite{Ner-Doug} for the complex analog of the Jacobi transform).
It is interesting to find a collection
of  evaluations of $J_{\alpha,\beta} f$ for some functions $f$. This section contains 
few examples.

\sm

{\it The transform $J_{\alpha,\beta}$ sends the function
	\begin{equation*}
	\left(\tfrac12 +i x \right)^{-p} 	\left(\tfrac12 -i x \right)^{-q}
	\end{equation*}
	to the function } 
	\begin{equation}
	2\pi \Gamma\left(p+q+\osigma-\tfrac12\right)\,\HH
	\left[\begin{matrix}
	\frac{1-\alpha-i\beta}{2}+\osigma+\ot,\,\,\frac{1+\alpha+i\beta}{2}+\osigma+\ot,\,\,1-q-\ot\\
	1-\frac{\alpha-i\beta}{2}+\ot,\,\,1+\frac{\alpha-i\beta}{2}+\ot,\,\,p+\frac12+\osigma+\ot
	\end{matrix};1
	\right].
	\label{eq:3H3}
	\end{equation}
	
To verify this, we must evaluate the integral
\begin{multline*}
\int_{-\infty}^{\infty} \left(\tfrac12 +i x \right)^{-p} 	\left(\tfrac12 -i x \right)^{-q}
\times\\\times
\ov {
 \left(\tfrac12+ix\right)^{t} \left(\tfrac12-ix\right)^{-1/2-t-\osigma}
 \!\!
\H\left[\begin{matrix}
\frac{ 1-\alpha+i \beta}2+\osigma+\ot,\,\,\frac{1+\alpha-i \beta}2+\osigma+t\\
1-\frac{\alpha+i\beta}{2}+\ot,\,\,1+\frac{\alpha+i\beta}{2}+\ot
\end{matrix};-\frac{\frac12+ix}{\frac12-ix}\right]}\,dx.
\end{multline*}	
We expand $\H$ into a series. Integrating term-wise with the formula
$$
\int_{-\infty}^\infty\frac{dx}{\left(\frac12+ix\right)^\mu \left(\frac12-ix\right)^\nu}
=\frac{2\pi\,\Gamma(\mu+\nu-1)}{\Gamma(\mu)\,\Gamma(\nu)},
$$
we come to \eqref{eq:3H3}.

\sm

{\sc Remark.} The functions (\ref{eq:3H3}) are bilateral version of Hahn functions, which were considered
in \cite{GKR}.
\hfill $\boxtimes$

\sm

Next, for cases $q=\pm \frac{\alpha-i\beta}{2}$ and for $p=\pm \frac{\alpha+i\beta}{2}$
 the expression (\ref{eq:3H3})
can be simplified. For instance, {\it the transformation
$J_{\alpha,\beta}$ sends a function
$$
\left(\tfrac12 +i x \right)^{-p} 	\left(\tfrac12 -i x \right)^{\alpha/2-i\beta/2}
$$
to}
\begin{equation*}
\frac{-2\sin\pi(\frac{\alpha-i\beta}{2}+\ot)}{\Gamma(p+\frac{\alpha+i\beta}{2}) \, 
	\Gamma(p-\frac{\alpha+i\beta}{2})}\cdot
\frac{\Gamma(p+\frac{-\alpha+i\beta-1}{2}+\osigma)\,\Gamma(p+\frac{-\alpha+i\beta-1}{2}+\osigma)}
{\Gamma(\frac12 -\alpha-\osigma)\,\Gamma(\frac12+i\beta-\osigma)}.
\end{equation*}

To establish this statement, we substitute $q=-\alpha/2+i\beta/2$ to \eqref{eq:3H3}. Then we get a function of the form
$$
\HH\left[\begin{matrix}
a_1,a_2,c\\b_1,b_2,c
\end{matrix};1\right]=\frac{\sin\pi c}{\pi}\H\left[\begin{matrix}
a_1,a_2\\b_1,b_2
\end{matrix};1\right]
$$
and apply the Dougall formula \eqref{eq:dougall}.

 In a similar way we set
$p=\alpha/2+i\beta/2$
and observe that {\it our transform sends
$$
\left(\tfrac12 +i x \right)^{-(\alpha+i\beta)/2} 	\left(\tfrac12 -i x \right)^{-q}
$$
to}
$$
\frac{2\cos\pi(\frac{\alpha+i\beta}{2}+\osigma+\ot)}
{\Gamma(q+\frac{\alpha-i\beta}{2})\,\Gamma(q-\frac{\alpha-i\beta}{2})}
\cdot\frac{\Gamma(q+\frac{\alpha+i\beta-1}{2}+\osigma)\Gamma(q+\frac{\alpha+i\beta-1}{2}-\osigma)}
{\Gamma(\frac12+\alpha-\osigma)\Gamma(\frac12+i\beta-\osigma)}.
$$

Two remaining cases are similar.

\tt
\noindent
Math. Dept., University of Vienna\& Pauli Institute \\
\&Institute for Theoretical and Experimental Physics (Moscow); \\
\&MechMath Dept., Moscow State University;\\
\&Institute for Information Transmission Problems;\\
URL: http://mat.univie.ac.at/$\sim$neretin/


\begin{thebibliography}{cc}
	
	\bibitem{Aga}
	Agarwal R. P. {\it General transformations of bilateral cognate trigonometrical series of ordinary hypergeometric type}. Canad. J. Math. 5 (1953), 544-553.

\bibitem{AAR}
Andrews G. E., Askey R., Roy R. {\it  Special functions.}
Cambridge University Press, Cambridge, 1999.	

\bibitem{Ask-Rom}	
 Askey R., {\it An integral of Ramanujan and orthogonal polynomials}, J. Indian Math. Soc. 51 
(1987), 27-36. 
	

	
	\bibitem{BS}
Berezin F.A., Shubin M. A. {\it The Schr\"odinger equation}.   Kluwer, Dordrecht, 1991.	

\bibitem{Che}
 Cherednik I., {\it Inverse Harish-Chandra transform and difference operators.} Internat. Math. Res. Notices 1997
(1997) 733--750.
	
\bibitem{vDE}	
	van Diejen J. F., Emsiz E. {\it Difference equation for the Heckman--Opdam hypergeometric function and its confluent Whittaker limit.} Adv. Math. 285 (2015), 1225-1240.
	
\bibitem{Dun}
Dunford N., Schwartz J. T. {\it Linear operators. Part II: Spectral theory. Self adjoint operators in Hilbert space.} John Wiley \& Sons, New York-London 1963.	
	
	\bibitem{HTF1}
	Erd\'elyi A., Magnus W., Oberhettinger F., Tricomi F. G.
	{\it Higher transcendental functions.} Vol.~I, 
	Based, in part, on notes left by Harry Bateman. McGraw-Hill, New York-Toronto-London, 1953. 
	
		\bibitem{HTF2}
	Erd\'elyi A., Magnus W., Oberhettinger F., Tricomi F. G.
	{\it Higher transcendental functions.} Vol.~II, 
	Based, in part, on notes left by Harry Bateman. McGraw-Hill, New York-Toronto-London, 1953. 
	
	\bibitem{FYa}
	Faddeev L. D., Yakubovski\u{\i} O. A.
	{\it Lectures on quantum mechanics for mathematics students.}
 Amer. Math. Soc., Providence, RI, 2009.
 
 \bibitem{Gro}
 Groenevelt W.
 {\it Wilson function transforms related to Racah coefficients.}
 Acta Appl. Math. 91 (2006), no. 2, 133-191.
	
	\bibitem{GKR}
	Groenevelt W., Koelink E., Rosengren H. {\it Continuous Hahn functions as Clebsch-Gordan coefficients,} in {\it Theory and applications of special functions,} 221-284, Dev. Math., 13, Springer, New York, 2005.
	
	\bibitem{Koe} Koekoek R.,  Swarttouw R. F. 
	{\it The Askey-scheme of hypergeometric orthogonal polynomials and its $q$-analogue.}
	Delft University of Technology,
	Faculty of Information Technology and Systems,
	Department of Technical Mathematics and Informatics,
	Report no. 98-17,	1998. Preprint {\tt https://arXiv.org/pdf/math/9602214.pdf}
	
	
	\bibitem{Koo1}
	 Koornwinder T. H., {\it Jacobi functions and analysis on noncompact symmetric spaces}, in  {\it Special functions: Group theoretical aspects and applications}, Reidel, Dordrecht 1984, pp. 1-85.

\bibitem{Koo2}
 Koornwinder T. H., {\it Special orthogonal polynomial systems mapping to each other by the Fourier-Jacobi transform}, Lecture Notes in Math., vol. 1171, Springer-Verlag, Berlin 1985, pp. 174-183.
 
 
 \bibitem{Mol}
	 Molchanov  V. F., {\it Tensor products of unitary representations of the three-dimensional Lorentz group}, Math. USSR-Izv., 15:1 (1980), 113-143.
	
	\bibitem{MN} Molchanov V. F., Neretin Yu. A.
{\it A pair of commuting hypergeometric operators on the complex plane and bispectrality.}
J. of Spectral Theory (2021), {\tt DOI: 10.4171/JST/349}.
	

\bibitem{Ner-index}
 Neretin Yu. A., {\it Index hypergeometric transform and imitation of analysis of Berezin kernels on hyperbolic spaces}, Sb. Math., 192:3 (2001), 403-432.

\bibitem{Ner-beta}
 Neretin Yu. A., {\it Beta-integrals and finite orthogonal systems of Wilson polynomials}, Sb. Math., 193:7 (2002), 1071-1089.

\bibitem{Ner-Jacobi}
Neretin Yu. A., {\it Some Continuous Analogs of the Expansion in Jacobi Polynomials and Vector-Valued Orthogonal Bases}, Funct. Anal. Appl., 39:2 (2005), 106--119.

\bibitem{Ner-add}
Neretin Yu. A., 
{\it Index hypergeometric integral transform.} Addendum to Russian translation
of Andrews G.E., Askey R., and Roy R., {\it Special Functions}, MCCME, 2013, 607--624;
English version: Preprint, {\tt arXiv:1208.3342}. 

\bibitem{Ner-imaginary}
Neretin, Yu. A. 
{\it Difference Sturm-Liouville problems in the imaginary direction.}
J. Spectr. Theory 3 (2013), no. 3, 237--269.

\bibitem{Ner-Doug}
Neretin, Yury A. {\it An analog of the Dougall formula and of the de Branges-Wilson integral.} Ramanujan J. 54 (2021), no. 1, 93-106.

\bibitem{Ole}
Olevski\u{\i} M. N.
{\it On the representation of an arbitrary function in the form of an integral with a kernel containing a hypergeometric function.} (Russian)
Doklady Akad. Nauk SSSR (N.S.) 69, (1949), 11--14.


\bibitem{Puck}
Pukanszky L. {\it On the Kronecker products of irreducible unitary representations of the
$2 \times 2$ real unimodular group.} Trans. Amer. Math. Soc., 100 (1961), 116-152.

\bibitem{Rom}
 Romanovski V., {\it Sur quelques classes nouvelles de polyn\'omes orthogonaux}, C. R. Acad. Sci., 
Paris, 188 (1929), 1023-1025. 


\bibitem{Sla}
Slater L. J. {\it Generalized hypergeometric functions.} Cambridge University Press, Cambridge, 1966.

\bibitem{Tit}
 Titchmarsh E. C., {\it Eigenfunction expansions with second-order differential operators}, Clarendon Press, Oxford, 1946.

\bibitem{Wey}
 Weyl H., {\it \"Uber gew\:onliche lineare Differentialgleichungen mit singul\"aren Stellen und ihre Eigenfunktionen (2 Note)}, Nachr. Konig. Gess. Wiss. G\"ottingen. Math.-Phys. (1910), 442-467; in H. Weyl, {\it Gesammelte Abhandlungen}, vol. 1, Springer-Verlag, Berlin 1968, pp. 222-247.
 
 \bibitem{Yak}
 Yakubovich, S. B. {\it Index transforms.}  World Scientific, River Edge, NJ, 1996.

\end{thebibliography}
\end{document}